\documentclass[12pt,leqno]{amsart}
\usepackage{amsmath,amsthm,a4}
\newtheoremstyle{fancy}{}{}{\itshape}{}{\textsc\bgroup}{.\egroup}{ }{}
\newcounter{theorem}
\theoremstyle{fancy}

\newtheorem{lem}[theorem]{Lemma}
\newtheorem{thm}[theorem]{Theorem}
\newcommand{\cref}[1]{Corollary~\ref{#1}}
\newcommand{\lref}[1]{Lemma~\ref{#1}}
\newcommand{\tref}[1]{Theorem~\ref{#1}}
\def\R{\mathbb R}       \def\Z{\mathbb{Z}}
\def\ve{\varepsilon} \def\la{\langle} \def\ra{\rangle}
\newcommand{\diam}{\operatorname{diam}}

\renewcommand{\Re}{\operatorname{Re}}
\newcommand{\Ric}{\operatorname{Ric}}

\begin{document}


\title{
Eigenvalues and Holonomy}
\author{
Werner Ballmann,
Jochen Br\"uning,
\and Gilles Carron}
\address{
Mathematisches Institut\\
Universit\"at Bonn\\
Beringstrasse 1, D-53115 Bonn\\}
\email{ballmann\@@math.uni-bonn.de}
\address{
Institut f\"ur Mathematik\\
Humboldt--Universit\"at\\
Rudower Chaussee 5, 12489 Berlin, Germany\\}
\email{bruening\@@mathematik.hu-berlin.de}
\address{
Departement de Mathematiques\\
Universite de Nantes\\
2 rue de la Houssiniere, BP 92208\\
44322 Nantes Cedex 03, France\\}
\email{Gilles.Carron@@math.univ-nantes.fr}

\thanks{The authors were partially supported
by SFB 256 (U Bonn) and SFB 288 (HU Berlin).}

\date{01.07.02}


\begin{abstract}
We estimate the eigenvalues of connection Laplacians
in terms of the non-triviality of the holonomy.
\end{abstract}

\maketitle


\section*{Introduction}

Let $S_L=\R/L\Z$ be a circle of length $L$
and $X$ be the oriented unit vector field on $S=S_L$.
Up to equivalence,
there is exactly one Hermitian line bundle, $E$, over $S$.
For a given complex number $z$ of modulus $1$,
there is, again up to equivalence, exactly one Hermitian
connection, $\nabla^E$, on $E$ with holonomy $z$ around $S$.

The Laplace operator $\Delta^E=(\nabla^E)^*\nabla^E$
is essentially self-adjoint as an operator in $L^2(E)$
with domain $C^2(E)$.
The spectrum of its closure is discrete
and consists of the eigenvalues
\begin{equation*}
      \frac{4\pi^2}{L^2}(\rho+k)^2 , \quad k\in\Z ,
\end{equation*}
where we write $z=\exp(2\pi i\rho)$.
The corresponding eigenspaces are spanned by the functions
$\exp(2\pi i(\rho+k)x/L)$.
We see that, for $z\ne1$, the spectrum does not contain $0$,
and that we can estimate the smallest eigenvalue
in terms of $L$ and $z$.

The aim of this paper is a correponding estimate for
Hermitian vector bundles over closed Riemannian manifolds
in higher dimensions.
The results of this paper are of importance in \cite{BBC},
but seem to be also of independent interest.

Let $M$ be a closed Riemannian manifold of dimension $n\ge2$.
Let $-(n-1)\kappa\le0$ be a lower bound for
the Ricci curvature of $M$, i.e. $\Ric_M\ge-(n-1)\kappa$,
and let $D$ be an upper bound for the diameter of $M$,
$\diam M\le D$.
Let $E\to M$ be a Hermitian vector bundle over $M$
and $\nabla^E$ be a Hermitian connection on $E$.
The kernel of the associated connection Laplacian
$\Delta^E=(\nabla^E)^*\nabla^E$
consists of globally parallel sections of $E$.
The estimates we obtain are in terms of quantitave
measures for the non-existence of parallel sections,
that is, in terms of the holonomy of $E$.

Assume first that $\nabla^E$ is flat
and that the holonomy of $\nabla$ is irreducible
(and nontrivial).
Recall that for each point $x\in M$,
the fundamental group $\pi_1(M,x)$ of $M$ at $x$
admits a {\em short basis}, that is,
a generating set represented by loops of length
at most $2\diam M$, see \cite{Gr}.
Hence for each point $x\in M$,
there is a constant $\alpha(x)>0$ such that
for all $v\in E_x$ there is a smooth unit speed
loop $c:[0,l]\to M$ at $x$ of length $l\le2\diam M$
with holonomy $H_c$ satisfying
\begin{equation*}
  | H_c(v) - v | \ge \alpha(x) |v| .
\end{equation*}
There is also a constant $\ve(x)>0$ such that a loop
at $x$ has length $>2\diam M+\ve(x)$ unless
it is homotopic to a loop at $x$ of length $\le2\diam M$.
It follows that for any point $y\in M$ of distance $<\ve/4$
to $x$,
the homotopy classes of loops of length $\le2\diam M$ at $y$
are represented by
concatenated curves of the form $c_{xy}^{-1}*c*c_{xy}$,
where $c_{xy}$ denotes a fixed minimal geodesic from $x$
to $y$ and $c$ is a loop at $x$ of length $\le2\diam M$.
Since $\nabla^E$ is flat, parallel translation does not
depend on homotopy classes.
It follows that for each point $y$ sufficiently close to $x$,
there is a loop $c$ of length $\le2\diam M$ at $y$
which has the same non-trivial holonomy as the loop
$c_{xy}*c*c_{xy}^{-1}$ at $x$.
In particular, we can choose the constants $\alpha(x)$
such that they have uniform positive lower bounds locally.
By the compactness of $M$,
there is a uniform constant $\alpha>0$ such that,
for all $x\in M$ and $v\in E_x$,
there is a smooth unit speed loop $c:[0,l]\to M$ at $x$
of length $l\le2\diam M$ with holonomy $H_c$ satisfying
\begin{equation}\label{e:holon}
  | H_c(v) - v | \ge \alpha |v| .
\end{equation}
Our first estimate is as follows.

\begin{thm}\label{t:estimflat}
Suppose that $\nabla^E$ is flat
and that the holonomy of $\nabla^E$ satisfies \eqref{e:holon}.
Then, for each eigenvalue $\lambda$ of $\Delta^E$,
\begin{equation*}
  \sqrt\lambda\,
  \exp\big(c_0\sqrt{\lambda+(n-1)\kappa}\,\diam M\bigr)
  \ge \frac{\alpha}{2\diam M}
\end{equation*}
with a constant $c_0=c_0(n,\sqrt\kappa D)$.
In particular,
\begin{equation*}
  \sqrt\lambda \ge \min \left\{ \frac{1}{c_0\diam M},
  \frac{\alpha}{2\diam M}
  \exp\left(-c_0\sqrt{(n-1)\kappa} \diam M -1 \right)
  \right\} .
\end{equation*}
\end{thm}

For each point $x\in M$ and unit vector $v\in E_x$,
let $\beta(v)$ be the supremum of the ratios $|H_c(v)-v|/L(c)$,
where the supremum is taken over all non-constant loops $c$
starting at $x$, $L(c)$ denotes the length of $c$,
and $H_c$ the holonomy along $c$.
Set
\begin{equation}\label{e:holon2}
  \beta := \inf \{ \beta(v) \mid v\in E, |v|=1 \} .
\end{equation}
Note that by the definition of the constant $\alpha$
in \eqref{e:holon}, we have $\beta\ge\alpha/2\diam M$.
In the general case, i.e. if $\nabla^E$ is not necessarily flat,
we have the following estimate.

\begin{thm}\label{t:estim}
There are positive constants $a=a(n)$
and $c_1=c_1(n,\sqrt\kappa D)$ such that,
for each eigenvalue $\lambda$ of $\Delta^E$,
\begin{equation*}
  \sqrt\lambda\, \exp\big(
  c_1\sqrt{\lambda+(n-1)\kappa+n^2r+n^2r^2/\beta^2}\,\diam M
  \bigr) \ge \frac{\beta}{a} ,
\end{equation*}
where $r$ is a uniform bound for the pointwise operator norm
of $R^E$.
In particular,
\begin{equation*}
  \sqrt\lambda \ge \min \left\{ \frac{1}{c_1\diam M}, 
  \frac{\beta}{a} \exp \left(
  -c_1\sqrt{(n-1)\kappa+n^2r+n^2r^2/\beta^2} \diam M - 1 \right)
  \right\} .
\end{equation*}
\end{thm}

\noindent
The constants $a$, $c_0$ and $c_1$ in Theorems \ref{t:estimflat}
and \ref{t:estim} can be determined explicitly.
Except for the factor $1/a$, \tref{t:estim} 
implies \tref{t:estimflat}.
On the other hand, the proof of \tref{t:estimflat} is more
elementary than the one of \tref{t:estim} and exposes the main
ideas more clearly.
Moreover, the constant $c_0$ is better than the constant $c_1$,
that is, $c_0\le c_1$.
Both proofs rely on a Sobolev inequality of Gallot \cite{Ga}
and Moser iteration.
In the proof of \tref{t:estim} we actually need an extension
of Moser's iteration technique.

If part of the holonomy is trivial, then the corresponding
space of parallel sections determines a subbundle $E'$ of $E$.
The above results then apply to the orthogonal complement $E''$
of $E'$ in $E$.
On the other hand, for a section $\sigma=\sum\phi_i\sigma_i$
in $E'$, where the sections $\sigma_i$ are parallel,
we have $\Delta^E\sigma=\sum(\Delta\phi_i)\sigma_i$,
and hence the usual eigenvalue estimates for the Laplace
operator on functions as for example in \cite{LY}
or \cite {Zh} apply.

\section*{An apriori estimate}

For the convenience of the reader and since we will need a
modification further on, we give a short account of Moser iteration
as applied in \cite{Li}, see also \cite{Ga,Au}.
It will give rise to the following infinite product,
\begin{equation}
  A(x,y,z) := \prod_{i=0}^{\infty}
  \left( x + \frac{yz^i}{\sqrt{2z^i-1}} \right)^{1/z^i} ,
\end{equation}
where $x,y>0$ and $z>1$. Note that
\begin{equation}
  A(tx,ty,z) = t^{z/(z-1)} A(x,y,z)
  \quad\text{for $t>0$}.
\end{equation}
We have $\sqrt{2z^i-1}>z^{i/2}$ and $\ln(1+yz^{i/2})<yz^{i/2}$,
by the assumptions on $y$ and $z$, hence
\begin{equation}\label{estimateA}
  A(1,y,z) \le \exp \big( \frac{y}{1-1/\sqrt z} \big) .
\end{equation}
This estimate will be good enough for the present purposes.
For more precise estimates,
see \cite[p.467]{Li} and \cite[p.7]{Au}.

Let $M$ be a closed Riemannian manifold of dimension $n$
and volume $V$.
Let $\nabla$ and $\Delta$ be the Levi-Civita connection
and the Laplace operator on functions of $M$, respectively.
Denote by $\Vert\cdot\Vert_p$ the $L^p$-norm with respect
to the {\em normalized} Riemannian measure of $M$.
Let $q>1$ be in the Sobolev range, i.e.
there are  positive constants $B$ and $C$ such
that $M$ satisfies the {\em Sobolev inequality}
\begin{equation}\label{sobolev}
  \Vert f \Vert_{2q} \le
  B \Vert f \Vert_2 + C V^{1/n} \Vert df \Vert_2 ,
\end{equation}
for all smooth functions $f$ on $M$.
Let $F\to M$ be a Hermitian vector bundle with a Hermitian
connection $\nabla^F$.
Let $\Delta^F$ be the associated connection Laplacian.

\begin{lem}\label{l:dgnm}
Let $\sigma\in L^2(M,F)$ be a smooth section.
Assume that (pointwise)
\begin{equation*}
  V^{2/n} \la \Delta^F\sigma,\sigma \ra
  \le \Lambda^2 |\sigma|^2
\end{equation*}
for some constant $\Lambda\ge0$. Then
\begin{equation*}
  \Vert \sigma \Vert_\infty
  \le
  A(B,C\Lambda,q) \, \Vert \sigma \Vert_{2} .
\end{equation*}
\end{lem}

\begin{proof}
We may assume that $V=1$.
The function $f_\ve=\sqrt{|\sigma|^2+\ve^2}$ is smooth,
and, by Kato's inequality and our assumption, respectively,
we have the pointwise estimate
\begin{equation}\label{e:m1}
  f_\ve \Delta f_\ve
  \le \Re \la \sigma,\Delta^F\sigma \ra
  \le \Lambda^2 |\sigma|^2
  \le \Lambda^2 f_\ve^2 .
\end{equation}
Let $k>1/2$. Since $df_\ve^k=kf_\ve^{k-1}df_\ve$,
\begin{equation}\label{e:m2}
\begin{split}
  \Vert df_\ve^k \Vert_2^2
  &= k^2 \la f_\ve^{k-1}df_\ve, f_\ve^{k-1}df_\ve \ra_2
   = \frac{k^2}{2k-1} \la df_\ve, df_\ve^{2k-1} \ra_2
  \\&=
  \frac{k^2}{2k-1} \la \Delta f_\ve, f_\ve^{2k-1} \ra_2
  \le \frac{\Lambda^2k^2}{2k-1} \int f_\ve^{2k}
  = \frac{\Lambda^2k^2}{2k-1} \Vert f_\ve \Vert_{2k}^{2k} .
\end{split}
\end{equation}
Using (\ref{sobolev}), we get
\begin{equation}\label{e:m3}
  \Vert f_\ve \Vert_{2kq}^k
  = \Vert f_\ve^k \Vert_{2q}
  \le B \Vert f_\ve \Vert_{2k}^{k}
  + \frac{C \Lambda k}{\sqrt{2k-1}} \Vert f_\ve \Vert_{2k}^{k} .
\end{equation}
By letting $\ve\to0$ we conclude
\begin{equation}\label{e:m4}
  \Vert \sigma \Vert_{2kq} \le
  \left(B + \frac{C\Lambda k}{\sqrt{2k-1}}\right)^{1/k}
  \Vert \sigma \Vert_{2k} .
\end{equation}
Iterating this inequality with $k=q^j$, $j=0,1,\ldots,$ we get
\begin{equation}\label{e:m5}
\begin{aligned}
  \Vert \sigma \Vert_{2q^{j+1}}
  &\le \left( B + \frac{C\Lambda q^j}{\sqrt{2q^j-1}}
    \right)^{1/q^j} \Vert \sigma \Vert_{2q^j} \\
  &\le \prod_{i=0}^{j}
  \left( B + \frac{C\Lambda q^i}{\sqrt{2q^i-1}}
    \right)^{1/q^i} \Vert \sigma \Vert_{2} .
\end{aligned}
\end{equation}
Now $\Vert\sigma\Vert_{2q^{j+1}} \to \Vert\sigma\Vert_\infty$
as $j$ tends to $\infty$.
Hence the lemma.
\end{proof}

\section*{Proof of \tref{t:estimflat}}

 From now on, we assume that $M$ satisfies
$\Ric_M\ge-(n-1)\kappa$ and $\diam M\le D$.
We will use the following Sobolev inequality.

\begin{lem}[Gallot \cite{Ga}]\label{l:gallot}
There is a positive constant $c'=c'(n,\sqrt\kappa D)$
such that, for all  $p\in[1,\frac{n}{n-1}]$
and all smooth functions $f$ on $M$,
\begin{equation*}
  \Vert f \Vert_{\frac{2p}{2-p}} \le
  \Vert f \Vert_2 + \frac{2c'}{2-p}\diam M \Vert df \Vert_2 .
\end{equation*}
\end{lem}

\noindent
In other words, if the assumptions of \lref{l:gallot} are
satisfied, then $M$ satisfies a Sobolev inequality of
the type (\ref{sobolev}) with
\begin{equation}\label{e:const}
  q = \frac{p}{2-p} , \quad B=1 , \quad\mbox{and}\quad
  C = \frac{2c'}{2-p}\diam M V^{-1/n} .
\end{equation}
Note also that the function $c'$ can be chosen to be equal to
\begin{equation}\label{e:gallot}
  c'(n,d) = \left\{ \frac1d \int_0^d \big(
  \frac12e^{(n-1)d}\cosh t + \frac1{nd}\sinh t \big)^{n-1} \,dt
  \right\}^{1/n}
\end{equation}
with $d=\sqrt\kappa D$, compare \cite{Ga}.

\begin{thm}\label{t:estimpar}
Suppose that $\nabla^ER^E=0$.
Then, for each eigenvalue $\lambda$ of $\Delta^E$,
\begin{equation*}
  \sqrt\lambda\,
  \exp\big(c_0\sqrt{\lambda+(n-1)\kappa+n^2r}\,\diam M\bigr)
  \ge \beta
\end{equation*}
with $r$ and $\beta$ as in \tref{t:estim}
and $c_0=c_0(n,\kappa\sqrt D)$.
\end{thm}

Recall that $\beta\ge\alpha/2\diam M$ and $r=0$
under the assumptions of \tref{t:estimflat}.
Hence \tref{t:estimpar} implies \tref{t:estimflat}.

\begin{proof}[Proof of \tref{t:estimpar}]
Let $\sigma$ be a nonzero section of $E$ with
$\Delta^E\sigma=\lambda\sigma$.
Let $x\in M$ and choose $\beta'<\beta$.
Then there is a unit speed loop $c:[0,l]\to M$ at $x$,
of length $l$, with holonomy $H_c:E_x\to E_x$ satisfying
\begin{equation*}
  | H_c(\sigma(x)) - \sigma(x) |
  \ge \beta' l | \sigma(x) | .
\end{equation*}
Let $F_1,\ldots,F_k:[0,l]\to E$ be a parallel orthonormal
frame along $c$.
Express $\sigma\circ c$ as a linear combination of this frame,
$\sigma\circ c=\sum\phi^i F_i$.
By the assumption on the holonomy, we have
\begin{equation*}
\begin{aligned}
  \beta' l | \sigma(x) | = \beta' l | \phi(0) |
  &\le | \phi(l) - \phi(0) |
   \le \int_0^l | \phi'| \,dt \\
  &\le \int_0^l | (\nabla^E\sigma) \circ c | \,dt
   \le l \Vert \nabla^E\sigma \Vert_\infty .
\end{aligned}
\end{equation*}
Since we use the normalized volume element for our norms,
this gives
\begin{equation}\label{e:1}
  \beta \Vert \sigma \Vert_2
  \le \beta \Vert \sigma \Vert_\infty
  \le \Vert \nabla^E\sigma \Vert_\infty .
\end{equation}
On the other hand, $\nabla^E\sigma$ is a one-form with
values in $E$, that is,
a section of the bundle $F=\Lambda^1(T^*M)\otimes E$.
This bundle inherits a connection, $\nabla^F$,
from the Levi--Civita connection $\nabla$ of $M$
and the connection $\nabla^E$ of $E$.
In terms of a local orthonormal frame $X_1,\ldots,X_n$ of $M$
and a further local vector field $Z$,
the corresponding Bochner formula is
\begin{multline}\label{e:bochner}
  (\Delta^F \nabla^E\sigma)(Z) = \\
  \nabla^E_Z(\Delta^E \sigma)
  - \nabla^E_{\Ric Z}\sigma
  - 2\sum R^E(X_i,Z)\nabla^E_{X_i}\sigma
  - \sum (\nabla^E_{X_i}R^E)(X_i,Z)\sigma ,
\end{multline}
see e.g. Lemma 3.3.1 of \cite{LR}.
In particular, since $\Delta^E\sigma=\lambda\sigma$
and $\nabla^ER^E=0$,
\begin{equation}\label{e:2}
  \la \Delta^F (\nabla^E\sigma), \nabla^E\sigma \ra \le
  \big( \lambda + (n-1)\kappa + 2n^2r \big) |\nabla^E\sigma|^2 ,
\end{equation}
where we are somewhat generous in the estimate of the
curvature term.
 From (\ref{e:2}) and Lemmas \ref{l:dgnm} and \ref{l:gallot},
where we choose $p=(n+2)/(n+1)$ and $q=(n+2)/n$
and $C$ as in \eqref{e:const}, we conclude that
\begin{equation*}
  \Vert \nabla^E\sigma \Vert_\infty
  \le A\big( 1,\frac{(2n+2)c'}{n}
  \sqrt{\lambda+(n-1)\kappa+2n^2r}\,\diam M,\frac{n+2}{n} \big)
  \, \Vert \nabla^E\sigma \Vert_2 .
\end{equation*}
Now
$\Vert\nabla^E\sigma\Vert_2=\sqrt{\lambda}\Vert\sigma\Vert_2$
since $\Delta^E\sigma=\lambda\sigma$.
In combination with (\ref{estimateA}) and (\ref{e:1}),
this proves the asserted inequality.
\end{proof}

\section*{Proof of \tref{t:estim}}

We cannot apply the previous argument directly
to prove \tref{t:estim}.
The reason is that, in general,
the Bochner formula \eqref{e:bochner} only gives the estimate
\begin{multline}\label{e:3}
  \la \nabla^E\sigma, \Delta^F\nabla^E\sigma \ra \le
  \big( \lambda + (n-1)\kappa + n^2r \big) |\nabla^E\sigma|^2 \\
  - \sum_{i,j} \big\langle (\nabla^E_{X_i}R^E)(X_i,X_j)\sigma
  + R^E(X_i,X_j)\nabla^E_{X_i}\sigma,\nabla^E_{X_j}\sigma \big\rangle .
\end{multline}
Note that we distributed the terms arising
from $2\sum R^E(X_i,Z)\nabla^E_{X_i}\sigma$ in \eqref{e:bochner}
to both terms on the right hand side in \eqref{e:3}.
Now \eqref{e:3} involves $\sigma$ on the right hand side.
To overcome this problem, we have to modify the argument
in the proof of \lref{l:dgnm}.
We replace $\sigma$ there by the section $\nabla^E\sigma$
under discussion here
and set $f_\ve:=\sqrt{|\nabla^E\sigma|^2+\ve^2}$.
Instead of \eqref{e:m1},
we now have the pointwise estimate
\begin{equation*}
\begin{split}
  f_\ve \Delta f_\ve
  \le \Re \la \nabla^E\sigma,\Delta^F\nabla^E\sigma \ra
  &\le \big( \lambda + (n-1)\kappa + n^2r \big) f_\ve^2 \\
  - \sum_{i,j} &\big\langle (\nabla^E_{X_i}R^E)(X_i,X_j)\sigma
  + R^E(X_i,X_j)\nabla^E_{X_i}\sigma,\nabla^E_{X_j}\sigma \big\rangle .
\end{split}
\end{equation*}
Let $k\ge 1$. Then
\begin{equation*}\label{e:mm4}
\begin{split}
  \int_M |df_\ve^k|^2 \le
  &\frac{k^2}{2k-1} \big( \lambda + (n-1)\kappa + n^2r \big)
    \int_M f_\ve^{2k} \\
  &- \frac{k^2}{2k-1} \int_M\sum_{i,j} \big\langle
  \nabla^E_{X_i}(R^E(X_i,X_j)\sigma),\nabla^E_{X_j}\sigma
  \big\rangle f_\ve^{2k-2} ,
\end{split}
\end{equation*}
where it is understood that we choose, for each point $x\in M$,
an orthonormal frame $X_1,\ldots,X_n$
with $(\nabla_{X_i}X_j)(x)=0$.
As in \cite{LR}, the divergence theorem gives
\begin{equation*}\label{e:mm5}
\begin{split}
  -\int_M &\sum_{i,j} \big\langle
  \nabla^E_{X_i}(R^E(X_i,X_j)\sigma), \nabla^E_{X_j}\sigma
  \big\rangle f_\ve^{2k-2} \\
  &= \int_M\sum_{i,j} \big\langle R^E(X_i,X_j)\sigma,
  \nabla^E_{X_i}\nabla^E_{X_j}\sigma \big\rangle f_\ve^{2k-2} \\
  &+ 2(k-1)\int_M f_\ve^{2k-3} \sum_{i,j}df_\ve(X_i)
  \big\langle R^E(X_i,X_j)\sigma,\nabla^E_{X_j}\sigma\big\rangle .
 \end{split}\end{equation*}
Now $R(X_i,X_j)=-R(X_j,X_i)$; therefore,
with the above choice of frames,
\begin{equation*}
  \sum_{i,j} \big\langle R^E(X_i,X_j)\sigma,
  \nabla^E_{X_i}\nabla^E_{X_j}\sigma \big\rangle
  = \frac{1}{2}\sum_{i,j}\big|R^E(X_i,X_j)\sigma\big|^2 .
\end{equation*}
Hence
\begin{equation*}\label{e:mm6}
\begin{split}
  -\int_M &\sum_{i,j} \big\langle
  \nabla^E_{X_i}(R^E(X_i,X_j)\sigma),\nabla^E_{X_j}\sigma
  \big\rangle f_\ve^{2k-2} \\
  &\le \frac{n^2r^2}{2}\int_M |\sigma|^2f_\ve^{2k-2}+2(k-1) nr
  \int_M |\sigma|f_\ve^{2k-2}|df_\ve| \\
  &\le \frac{n^2r^2}{2} \|\sigma\|_{\infty}^2\int_Mf_\ve^{2k-2}
  + 2\frac{k-1}{k}nr\| \sigma\|_{\infty}
    \int_Mf_\ve^{k-1}|df_\ve^k| \\
  &\le \frac{n^2r^2}{2} \|\sigma\|_{\infty}^2\int_Mf_\ve^{2k-2}
  + 2nr \| \sigma\|_{\infty} \int_Mf_\ve^{k-1}|df_\ve^k| .
\end{split}
\end{equation*}
But
\begin{multline*}\label{e:mm7}
  \frac{2k(k-1)}{2k-1}nr\| \sigma\|_{\infty}\int_Mf_\ve^{k-1}|df_\ve^k|
  \le \\
  \frac{1}{2}\int_M|df_\ve^k|^2
  +\left(\frac{k(k-1)}{2k-1}\right)^22n^2r^2\|
        \sigma\|_{\infty}^2\int_Mf_\ve^{2k-2}
\end{multline*}
and $\|\sigma\|_{\infty}\le\|\nabla^E\sigma\|_{\infty}/\beta
    \le\|f_\ve\|_{\infty}/\beta$ , hence
\begin{equation}
\begin{split}
  \Vert df_\ve^k \Vert_2^2
  &\le
  \frac{2k^2}{2k-1} \left( \lambda + (n-1)\kappa + n^2r
  + \left(\frac{1}{2}+\frac{2(k-1)^2}{2k-1}\right)
    \frac{n^2r^2}{\beta^2} \right)
    \|f_\ve\|_{\infty}^2 \Vert f_\ve^{k-1} \Vert_2^2 \\
  &\le
  2k^2 \left( \lambda + (n-1)\kappa + n^2r
  + \left(\frac{1}{2}+\frac{2(k-1)^2}{(2k-1)^2}\right)
    \frac{n^2r^2}{\beta^2} \right)
    \|f_\ve\|_{\infty}^2 \Vert f_\ve^{k-1} \Vert_2^2 .
\end{split}
\end{equation}
Set
\begin{equation*}
  L^2:=2\big( \lambda + (n-1)\kappa + n^2r+n^2r^2/\beta^2 \big) .
\end{equation*}
Since $k\ge1$, we now have, instead of \eqref{e:m2},
\begin{equation*}
  \Vert df_\ve^k \Vert_2^2
  \le L^2k^2\Vert f_\ve \Vert_\infty^2
  \Vert f_\ve\Vert_{2k-2}^{2k-2}
\le L^2k^2\Vert f_\ve \Vert_\infty^2
  \Vert f_\ve \Vert_{2k}^{2k-2} .
\end{equation*}
Using \lref{l:gallot} with $q=(n+2)/n$
and $C$ as before, we replace \eqref{e:m3} by
\begin{equation*}
\begin{split}
  \Vert f_\ve \Vert_{2kq}^k
  = \Vert f_\ve^k \Vert_{2q}
  &\le \Vert f_\ve \Vert_{2k}^{k} + C L k
  \Vert f_\ve \Vert_\infty
  \Vert f_\ve \Vert_{2k}^{k-1} \\
  &\le (1 + C L k) \Vert f_\ve \Vert_\infty
  \Vert f_\ve \Vert_{2k}^{k-1} .
  \end{split}
\end{equation*}
Instead of \eqref{e:m3}, we conclude now, by letting $\ve\to0$,
that
\begin{equation}\label{e:mm3}
  \Vert \nabla^E\sigma \Vert_{2kq} \le
  \left(1 + C L k\right)^{1/k}
  \Vert \nabla^E\sigma \Vert_\infty^{1/k}
  \Vert \nabla^E\sigma \Vert_{2k}^{1-1/k} .
\end{equation}
As in \eqref{e:m4}, we iterate this inequality
with $k=q^j$, but now only for $j=1,2\ldots$
since \eqref{e:mm3} is useless in the case $k=1$.
Setting $p_i:=1-1/q^i$, we get
\begin{equation*}
\begin{aligned}
  \Vert \nabla^E\sigma \Vert_{2q^{j+1}}
  &\le
  \left( 1 + C L q^j \right)^{1/q^j}
    \Vert \nabla^E\sigma \Vert_\infty^{1-p_j}
    \Vert \nabla^E\sigma \Vert_{2q^j}^{p_j} \\
  &\le \prod_{i=1}^{j}
  \left( 1 + C L q^i
  \right)^{p_{i+1}\cdot\ldots\cdot p_j/q^i}
  \Vert \nabla^E\sigma \Vert_\infty^{1-p_1\cdot\ldots\cdot p_j}
  \Vert \nabla^E\sigma \Vert_{2q}^{p_1\cdot\ldots\cdot p_j} \\
  &\le \prod_{i=1}^{j}
  \left( 1 + C L q^i
  \right)^{1/q^i}
  \Vert \nabla^E\sigma \Vert_\infty^{1-p_1\cdot\ldots\cdot p_j}
  \Vert \nabla^E\sigma \Vert_{2q}^{p_1\cdot\ldots\cdot p_j} ,
\end{aligned}
\end{equation*}
where we use, for the latter inequality, that  $0<p_i<1$ and
that $x^p\le x$ if $x\ge1$ and $0<p<1$.
The limit
\begin{equation}\label{e:epsil}
  \ve = \ve(n):= \prod_{i=1}^{\infty} p_i
\end{equation}
exists and satisfies $0<\ve<1$. Moreover, using the inequality
\begin{equation*}
  1+C L q^i \le (1+C L)q^{i}
\end{equation*}
we obtain
\begin{equation*}
  \prod_{i=1}^{\infty}
  \left( 1 + C L q^i \right)^{1/q^i}
  \le (1 + C L)^{\sum_{i=1}^\infty 1/q^i} 
	\cdot q^{\sum_{i=1}^\infty i/q^i}
  \le a_1(n) e^{b(n) CL} 
\end{equation*}
with $a_1(n)=q^{\sum_{i=1}^\infty i/q^i}$ 
and $b(n)=\sum_{i=1}^\infty 1/q^i$.
We conclude that
\begin{equation*}
  \Vert \nabla^E\sigma \Vert_\infty \le a_2(n) 
  \exp \big(b(n)C L/\ve(n)\big) \Vert\nabla^E\sigma\Vert_{2q} 
\end{equation*}
with $a_2(n)=a_1(n)^{1/\ve(n)}$.
We also have
\begin{align*}
  \Vert \nabla^E\sigma \Vert_{2q}
  &\le \Vert \nabla^E\sigma \Vert_{2}^{1/q}
    \cdot \Vert \nabla^E\sigma \Vert_{\infty}^{(q-1)/q} \\
  &\le \Vert \nabla^E\sigma \Vert_2^{n/(n+2)}
    \cdot \Vert \nabla^E\sigma \Vert_{\infty}^{2/(n+2)} ,
\end{align*}
where we recall that $q=(n+2)/n$. Hence finally
\begin{equation*}
  \Vert \nabla^E\sigma \Vert_\infty \le a(n)
  \exp\big((n+2)b(n)C L/(n\ve(n))\big) \Vert\nabla^E\sigma\Vert_2 
\end{equation*}
with $a(n)=a_2(n)^{(n+2)/n}$.
The rest of the argument is as before.



\end{document}